\documentclass[11pt,epsfig]{article}
\usepackage{amsmath}
\usepackage{amsthm}
\usepackage{epsfig}
\usepackage{amsfonts}
\usepackage[margin=1.1in]{geometry}

\newtheoremstyle{thmm}{1.5ex plus 1ex minus
.2ex}{1.5ex plus 1ex minus
.2ex}{\rmfamily}{}{\bfseries}{}{1em}{}
\theoremstyle{thmm}
\newtheorem{theorem}{Theorem}[section]
\newtheorem{lemma}{Lemma}[section]

\newtheorem{corollary}{Corollary}[section]

\newcommand{\nn}{\nonumber}
\def \endproof{\vrule height8pt width 5pt depth
0pt}
\def\refe#1{(\ref{#1})}

\def\R{\mathbb{R}}

\begin{document}
\title{\bf Unconditionally optimal error analysis
of fully discrete Galerkin methods for general
nonlinear
parabolic equations}
\author{Buyang
Li\,\setcounter{footnote}{0}\footnote{Department
of Mathematics, City University of Hong Kong,
Kowloon, Hong Kong.
The work of the
authors was supported in part by a grant
from the Research Grants
Council of the Hong Kong Special Administrative
Region, China
(Project No. CityU 102005)\newline
\indent~ Email address:
{\tt libuyang@gmail.com},~  {\tt
maweiw@math.cityu.edu.hk} }
~~and ~Weiwei Sun\footnotemark[1]}
\date{}
\maketitle

\begin{abstract}
The paper focuses on unconditionally optimal error
analysis
of the fully discrete Galerkin finite element
methods for a general nonlinear
parabolic system in $\R^d$ with $d=2,3$.
In terms of a corresponding time-discrete system of
PDEs as proposed in \cite{LS1}, we
split the error function into two parts,
one from the temporal discretization and one the spatial discretization.
We prove that the latter is $\tau$-independent and the
numerical solution is bounded
in the $L^{\infty}$ and $W^{1,\infty}$ norms by the inverse inequalities.
With the boundedness of the numerical solution, optimal error estimates can
be obtained unconditionally in a routine way.
Several numerical examples in two and three
dimensional spaces are given to
support our theoretical analysis.
\end{abstract}

{\bf Key words:} Optimal error estimates,
unconditional stability,
Galerkin, nonlinear parabolic system
\medskip


\section{Introduction}
\setcounter{equation}{0}
There are several numerical approximations schemes in the
time direction for the numerical solution of
nonlinear parabolic equations (systems).
Linearized (semi)-implicit schemes are the most
popular ones
since, at each time step, the schemes only require
the solution of a linear system.
However, time-step size restriction condition
is always a key issue
in analysis and computation.  For many nonlinear
parabolic systems, error
analysis of finite element methods (or finite
difference method) with linearized
semi-implicit schemes in the time direction often
requires certain time-step conditions. See \cite{AG,He,KL,Liu1, shen} for the Navier-Stokes
equations,
\cite{EL,Zhao} for the nonlinear Joule heating
problems,
\cite{DEW,Duran,EW,HLS,Wang} for flows in porous media,
\cite{CL,EM, WHS} for viscoelastic fluid flow,
\cite{MS,ZM} for the KdV equations,
\cite{CH,MH} for the Ginzburg-Landau equations,
\cite{ADK,BC,Tou,Zou} for the nonlinear
Schr\"{o}dinger equations
and \cite{DM,FHL,WML} for some other equations.
Such time-step size restrictions may result in the use
of an unnecessarily small time-step size and extremely time-consuming in practical
computations.
To study the error estimate of linearized
(semi)-implicit schemes,
the boundedness of numerical solution (or error
function)
in the $L^{\infty}$ norm or a stronger norm is often
required.
If a priori estimate for the numerical solution in
such a norm cannot be provided, one
may employ the induction method with an inverse
inequality to bound
the numerical solution, such as
\begin{equation}
\| U_h^n - R_hu^n \|_{L^\infty} \leq
C h^{-d/2} \| U_h^n - R_hu^n \|_{L^2} \leq
C h^{-d/2} (\tau^{m} + h^{r+1}) ,
\end{equation}
where $U_h^n$ is the finite element solution, $u$
is the exact solution and $R_h$ is
certain projection operator. A time-step size
restriction arises immediately from the above
inequality,
particularly for problems in three
dimensional spaces.
Most previous works follow this idea.  A new
approach for unconditionally optimal error
analysis of a linearized Galerkin FEM was pesented
in our recent work  \cite{LS1}, also see
\cite{Li}, where
the error function is split into two parts,
the spatially discrete error and the temporally
discrete error,
\begin{equation}
\| u^n -U_h^n \| \le \| u^n -U^n \| + \| U^n -
U_h^n \|,
\end{equation}
where $U^n$ is the solution of a corresponding
time-discrete parabolic equations (or elliptic
equations).
Optimal estimates for the second term can be
obtained unconditionally in a traditional way if
suitable regularity of the solution of the
time-discrete system can be proved.
More recently, unconditionally optimal error
estimates
were established for a nonlinear equation from
incompressible miscible flow
in porous media.  In \cite{LS1,LS2}, analysis
waw given only for a linear FEM and
a low-order Galerkin-mixed FEM, respectively.

In this paper, we consider a general nonlinear
parabolic equation (or system)
\begin{align}
&\frac{\partial u}{\partial
t}-\nabla\cdot(\sigma(u)\nabla u)=g(u,\nabla
u,x,t)
\label{e-parab-1}
\end{align}
in a bounded and smooth domain $\Omega$ in $\R^d$
($d=2$ or $3$) with the
boundary condition
\begin{align}
\label{BC}
\begin{array}{ll}
u=0~~
&\mbox{on}~~\partial\Omega
\end{array}
\end{align}
and the initial condition
\begin{align}
\label{IniC}
\begin{array}{ll}
u(x,0)=u_0(x)~~
&\mbox{for}~~x\in\Omega ,
\end{array}
\end{align}
where $g\in C^2(\R)$ is a general nonlinear source.
The general equation is of stronger nonlinearity
than those in \cite{LS1, LS2} and many physical
equations are included. We apply linearized
backward Euler Galerkin method with $r$-order
finite element
approximation ($r\ge 1$) for the general nonlinear
system.
We focus our attention  on the unconditional
convergence (stability) and optimal error
estimates of
the linearized Galerkin FEMs.   A key to our
analysis is
the a priori estimate of the numerical solution.
We apply the splitting technique proposed in
\cite{LS1,LS2} to bound the
numerical solution $U_h^n$ in $L^{\infty}$-norm
and $W^{1,\infty}$-norm, such as
\begin{align}
\| U_h^n \|_{W^{1,\infty}}
& \le \| R_h U^n \|_{W^{1,\infty}} +
\| U_h^n - R_h U^n \|_{W^{1,\infty}} \nn \\
&\le \| R_h U^n \|_{W^{1,\infty}} + C h^{-d/2} \|
U_h^n - R_h U^n \|_{H^1}  \nn \\
& \le C + C h^{-d/2}
h^{k}
\end{align}
where $k>d/2$. Then with the boundedness of $\| U_h^n \|_{W^{1,\infty}}$,
optimal error estimates can be easily established
unconditionally
in the routine way of FEM error analysis.

The paper is organized as follows.
In Section 2, we present linearized backward Euler
Galerkin FEMs for the general nonlinear
parabolic equations \refe{e-parab-1}-\refe{IniC}
and introduce our notations.
In Section 3, we prove the boundedness of the
numerical solution in the
$W^{1,\infty}$ norm in terms of a corresponding time-discrete system, in which a rigorous
analysis on the regularity of the solution to the
time-discrete PDEs is given.
Due to the boundedness of the numerical solution in the
$W^{1,\infty}$ norm, we present unconditionally optimal error estimates
in Section 4 in a simple and routine way.
Numerical examples in two and
three-dimensional spaces are presented in
Section 5. Numerical results confirm our
theoretical analysis and show that no time-step condition is needed.

\section{Fully discrete Galerkin FEMs}
\setcounter{equation}{0}
Let $\pi_h$ be a regular division of $\Omega$ into
triangles
$T_j$, $j=1,\cdots,M$, in $\R^2$ or tetrahedras in
$\R^3$, and let
$h=\max_{1\leq j\leq M}\{\mbox{diam}\,T_j\}$
denote the mesh size.
For a triangle $T_j$ with two corners (or a
tetrahedra with three corners) on the
boundary, we let $\widetilde T_j$ denote the
triangle with one curved
side (or a tetrahedra with one curved face).  For
an interior triangle, we
simply set $\widetilde T_j$ as $T_j$ itself.
Finite element spaces on $\{ \widetilde T_j \}$
have been well defined,  $e.g.$, see \cite{Tho,
Zlamal}. For
a given triangular (or tetrahedral) division of
$\Omega$, we define the finite element space
\begin{align*}
&\widehat V_{h}=\{v_h\in C(\overline\Omega_h):
v_h|_{T_j}
\mbox{~is~a~polynomial of degree $r$ and $v_h=0$
on $\partial\Omega_h$}\}
\end{align*}
so that $\widehat V_h$ is a subspace of
$H^1_0(\Omega_h)$.
Let $G:\Omega_h\rightarrow \Omega$ be a coordinate
transformation
such that both $G$ and $G^{-1}$ are Lipschitz
continuous and,
for each triangle $T_j$, $G$ maps $T_j$ one-to-one
onto $\widetilde T_j$ \cite{Zlamal}.
We define an operator ${\cal
G}:L^2(\Omega_h)\rightarrow L^2(\Omega)$ by
${\cal G} v(x)=v(G^{-1}(x))$ for $x\in\Omega$.
Then we set
\begin{align*}
&V_{h}=\{{\cal G}v_h: v_h\in \widehat V_{h}\} .
\end{align*}
Easy to see that $V_h$ is a finite element
subspace of $H^1_0(\Omega)$ and
\begin{align*}
&\|w-\Pi_hw\|_{L^p}\leq
C\|w\|_{W^{r+1,p}(\Omega)}h^{r+1},
\quad\mbox{for}~r\geq 1 ,~2\leq p\leq\infty,
\end{align*}
where $\Pi_h={\cal G}\widehat\Pi_h{\cal G}^{-1}$
and
$\widehat\Pi_h:C_0(\overline\Omega_h)\rightarrow
\widehat V_h$ is the
Lagrangian interpolation operator of degree $r$.

Let $0=t_0<t_1<\cdots<t_N=T$ be a
uniform partition of the time interval $[0,T]$
with $t_n=n\tau$ and
let $u^n = u(x,t_n)$ for $n=0,1,\cdots,N$.
For a sequence of functions $\{ f^n \}_{n=0}^N$,
we define
\begin{align}
& D_\tau f^{n+1}=\frac{f^{n+1}-f^n}{\tau},
\quad\mbox{for $n=0,1,\cdots,N-1$.}
\label{dfnsjlwejriow}
\end{align}
A simple linearized backward Euler Galerkin method
for the problem
(\ref{e-parab-1})-(\ref{IniC})
is to seek $U^{n+1}_h\in V_h$, $0\leq n\leq N-1$,
such that
\begin{align}
&\big(D_\tau
U^{n+1}_h,v\big)+\big(\sigma(U^n_h)\nabla
U^{n+1}_h,\nabla v\big)
=\big(g(U^n_h,\nabla U^n_h,x,t),v\big)
\label{FDFEM}
\end{align}
for any $v\in V_h$, with the initial condition
$U^0_h= \Pi_hu_0$ for $r\geq 2$ and
$U^0_h=R_h^1u_0$ for $r=1$,
where $R_h^1$ is a projection operator defined in
Section 3.2.

With a linear approximation to the nonlinear
source term, an alternative linearized scheme
is defined by
\begin{align*}
&\big(D_\tau
U^{n+1}_h,v\big)+\big(\sigma(U^n_h)\nabla
U^{n+1}_h,\nabla v\big)
=\big(g_0^n,v\big ) + \big(g_1^n U^{n+1}_h, v\big
)+  \big(g_2^n \cdot \nabla U^{n+1}_h, v\big
),\quad\forall~v\in V_h,
\end{align*}
where $g_0^n = g(U^n_h,\nabla U^n_h,x,t)$,
$g_1^n = \nabla_1 g(U^n_h,\nabla U^n_h,x,t)$ and
$g_2^n= \nabla_2 g(U^n_h,\nabla U^n_h,x,t)$, with
$\nabla_1 g$ and $\nabla_2 g$
denoting the gradient of $g$ with respect to the
components
$U$ and $\nabla U$, respectively. The
corresponding linearized Crank-Nicolson schemes
can be defined similarly with classical
extrapolations \cite{DFJ}.

In this paper, we only focus our attention on the
linearized scheme \refe{FDFEM}.
The analysis presented in this paper can be
extended to the second linearized scheme
and many other schemes.
We assume that $g\in
C^2(\R\times\R^d\times\overline\Omega\times[0,T])$
and $\sigma\in C^2(\R)$ satisfies the weak
ellipticity condition
\begin{align}
\label{Strongellip}
\sigma(s)>0\quad \mbox{for}~~s\in\R .
\end{align}

\section{Boundedness of the numerical solution}
\setcounter{equation}{0}
In this section,
we assume that the solution to the problem
(\ref{e-parab-1})-(\ref{BC}) exists and satisfies
that
\begin{align}
\label{regularity}
&\|u\|_{L^\infty((0,T);H^3)}
+\|\partial_t
u\|_{L^\infty((0,T);L^2)}+\|\partial_t
u\|_{L^2((0,T);H^2)}
+\|\partial_{tt}u\|_{L^2((0,T);L^2)}  \leq M ,
\end{align}
for some positive constant $M$, and
we prove the following theorem.

\begin{theorem}\label{ErrestFEMSol}
{\it
Suppose that the system
\refe{e-parab-1}-\refe{IniC} has a unique solution
$u$ satisfying
the regularity condition \refe{regularity}. Then
there exist positive constants $C$ and $h_0$,
independent of $n$ and $h$,
such that the finite element system (\ref{FDFEM})
admits a unique solution
$\{ U^n_h \}_{n=1}^N$ when $h<h_0$, and
\begin{align}\label{optimalL2est}
\|U_h^n\|_{L^\infty}+\|\nabla
U_h^n\|_{L^\infty}\leq C .
\end{align}
}
\end{theorem}
\medskip

To prove Theorem \ref{ErrestFEMSol}, we introduce
a corresponding time-discrete
equation as proposed in \cite{LS1,LS2}:
\begin{align}
&D_\tau U^{n+1}-\nabla\cdot(\sigma(U^n)\nabla
U^{n+1})=g(U^n,\nabla U^n,x,t_n)
\label{e-TDparab-1},
\end{align}
with the boundary condition $U^{n+1}=0$ on
$\partial\Omega$ and the initial
condition $U^0=u_0$.

In the following two subsections, we estimate the
error functions
$u^n-U^n$ and $U^n- U_h^n$, respectively, where
$U^n$ is the solution of the time-discrete system
\refe{e-TDparab-1}.

For the simplicity of notations, we denote by $C$
a generic positive
constant and by $\epsilon$ a generic small
positive constant,
which depend solely upon $M$, $\Omega$, $T$,
$\sigma$ and $g$,
and independent of $\tau$, $h$ and $n$.

\subsection{The time-discrete solution}
By the regularity assumption (\ref{regularity}),
we have
$u\in W^{1,\infty}$. We set
\begin{align*}
& K=\|u\|_{L^\infty(\Omega \times [0,T])}+\|\nabla
u\|_{L^\infty(\Omega \times [0,T])}+2 \, .
\nn \\
& Q_K = [-K, K]^{d+1} \times \Omega \times [0,T]
\, .
\end{align*}
Then, by the regularity assumptions on $g$ and $\sigma$ and the ellipticity condition (\ref{Strongellip}), there exist positive constants $\sigma_K$ and
$C_K$ such that
for $|s|\leq K$ and $(\alpha, \beta, x, t) \in
Q_K$,
\begin{align}
&\sigma_K \leq \sigma(s)\leq C_K , \nn \\[5pt]
&|\sigma'(s)|+|\sigma''(s)|\leq C_K ,
\label{g-sigma}\\[-5pt]
&|g(\alpha,\beta,x,t)|
+|\partial_{\alpha} g(\alpha, \beta,x,t)|
+\sum_{j=1}^d|\partial_{\beta_j}g(\alpha,\beta,x,t)|
+\sum_{j=1}^d|\partial_{x_j}g(\alpha,\beta,x,t)|
\leq C_K \nn .
\end{align}

\begin{lemma}\label{h23}($H^{l}$-estimate of
elliptic equations \cite{ChenYZ})
{\it Suppose that $v$ is a solution of the
boundary value problem
\begin{align*}
& \Delta v = f, \quad \mbox{in}~~\Omega,
\\
& v = 0 , \quad \mbox{on}~~ \partial \Omega,
\end{align*}
where $\Omega \in \R^d$, $d=2,3$, is a smooth and
bounded domain. Then
\begin{equation}
\| v \|_{H^l} \le C \| f \|_{H^{l-2}} , \quad
l=2,3 \, .
\end{equation}
}
\end{lemma}

In this subsection, we explore the regularity of
the solution to
the time-discrete system (\ref{e-TDparab-1}) and
present an
error estimate for $u^n-U^n$.

\begin{theorem}\label{ErrestDisSol}
{\it Suppose that the system
\refe{e-parab-1}-\refe{IniC} has a unique
solution $u$ satisfying \refe{regularity}. Then
the time-discrete system (\ref{e-TDparab-1})
admits a
unique solution $\{ U^{n} \}_{n=0}^N$ such that
\begin{align}
& \max_{0 \le n \le N} \|U^n\|_{H^3}
+ \sum_{n=1}^N\tau\| D_\tau U^{n}\|_{H^2}^2
\leq C_0,
\label{ErrestDisSol211}
\end{align}
and
\begin{align}\label{TDErrEstLemm}
\begin{array}{ll}
&\displaystyle\max_{0\leq n\leq N}\|
u^{n}-U^n\|_{H^1} \leq C_0\tau
\end{array}
\end{align}
where $C_0$ is a positive constant independent of
$n$, $h$ and $\tau$.
}
\end{theorem}

\noindent{\it Proof}~~~
For the given $U^n$,  \refe{e-TDparab-1} can be
viewed as
a linear elliptic boundary value problem.
With the first inequality in (\ref{g-sigma}) and
classical theory of elliptic PDEs,
the equation \refe{e-TDparab-1} admits a unique
solution $U^{n+1}$
in $H^1$. Let $e^n=u^n-U^n$.  Here we only prove
the estimates
\refe{ErrestDisSol211}-\refe{TDErrEstLemm}.

First, we prove by mathematical induction the
inequality
\begin{align}
\label{fsdioo123}
\|U^n\|_{L^\infty}+\|\nabla
U^n\|_{L^\infty}<K,\quad\mbox{for}~~n=0,1,\cdots,N
\end{align}
under the condition $\tau<\tau_0$ for some
positive constant $\tau_0$.
Since $U^0=u_0$, the inequality \refe{fsdioo123}
holds for $n=0$.
Now we assume that the inequality holds for $0\leq
n\leq k$.

Let $e^n=u^n-U^n$. From
\refe{e-parab-1}-\refe{IniC} and
\refe{e-TDparab-1}, we see that
$e^{n+1}$ satisfies the equation
\begin{align}
&D_\tau e^{n+1}-\nabla\cdot(\sigma(U^n)\nabla
e^{n+1}) \label{err-TDparab-01}\\
&=R^{n+1}+\nabla\cdot([\sigma(u^n)-\sigma(U^n)]\nabla
u^{n+1})
+g(u^n,\nabla u^n,x,t)-g(U^n,\nabla U^n,x,t), \nn
\end{align}
with the boundary condition $e^{n+1}=0$ on
$\partial\Omega$
and the initial condition $e^0=0$,
where
\begin{align*}
R^{n+1}=\partial_tu^{n+1}
-D_\tau u^{n+1}+\nabla\cdot
[(\sigma(u^n)-\sigma(u^{n+1}))\nabla u^{n+1}]\\
+g(u^n,\nabla u^n,x,t)
-g(u^{n+1},\nabla u^{n+1},x,t)
\end{align*}
is the truncation error due to the time
discretization. By the regularity assumption
(\ref{regularity}), we have
\begin{align}\label{trdfsjhklh}
\max_{1\leq n\leq N}\|R^{n}\|_{L^2}\leq C,\qquad
\sum_{n=1}^{N}\tau
\|R^{n}\|_{L^2}^2\leq C\tau^2 .
\end{align}
Multiplying the equation (\ref{err-TDparab-01})
by
$e^{n+1}$, we obtain
\begin{align*}
D_\tau\biggl(\frac{1}{2}\|e^{n+1}\|_{L^2}^2\biggl)
&+ \frac{\sigma_K}{2} \| \nabla e^{n+1}
\|_{L^2}^2
\leq \biggl((\sigma(U^n)-\sigma(u^n))\nabla
u^{n+1}~,~\nabla e^{n+1}\biggl)
\\
& +\biggl(g(u^n,\nabla u^n,x,t)-g(U^n,\nabla
U^n,x,t), e^{n+1} \biggl)
+ (R^{n+1}, \, e^{n+1} ) \, .
\end{align*}
By \refe{g-sigma}, we have further
\begin{align*}
& |g(u^n,\nabla u^n,x,t)-g(U^n,\nabla U^n,x,t)|
\le C_K (|e^n|+|\nabla e^n|) ,
\\
&| \sigma(U^n) - \sigma(u^n) | \le C_K | e^n|
\, .
\end{align*}
It follows that
\begin{align*}
&
D_\tau\biggl(\frac{1}{2}\|e^{n+1}\|_{L^2}^2\biggl)+
\sigma_K \| \nabla e^{n+1} \|_{L^2}^2
\\
&\leq \epsilon \|\nabla e^{n+1} \|_{L^2}^2
+ C\epsilon^{-1}\|e^{n+1}\|_{L^2}^2 + \epsilon
\|\nabla e^{n} \|_{L^2}^2
+  C\epsilon^{-1}\|e^{n}\|_{L^2}^2
+C\epsilon^{-1}\|R^{n+1}\|_{L^2}^2  .
\end{align*}
By choosing $\epsilon<\sigma_K/4$ and using
Gronwall's inequality, there exists
$\tau_1 >0$ such that when $\tau \le \tau_1$
\begin{align}\label{dsfjkl}
&\| e^{n+1}\|_{L^2}\leq
C\biggl(\sum_{m=0}^{n}\tau
\|R^{m+1}\|_{L^2}^2\biggl)^{\frac{1}{2}}\leq C\tau
, \quad 0\leq n\leq k,
\end{align}
which implies that
\begin{align}\label{ion2390}
&\|U^{n+1}\|_{L^2}\leq \| u^{n+1}\|_{L^2}+\|
e^{n+1}\|_{L^2}  \leq C,\\
&\|D_\tau U^{n+1}\|_{L^2}\leq \|D_\tau
u^{n+1}\|_{L^2}+\|D_\tau e^{n+1}\|_{L^2}  \leq C
\, .
\end{align}
Applying Lemma \ref{h23} for the linear elliptic
equation
\refe{e-TDparab-1} with the induction assumption
gives the $H^2$ estimate
\begin{align}
 \| U^{n+1} \|_{H^2}
& \le
C \|D_\tau U^{n+1}\|_{L^2} + C \| \nabla
\sigma(U^n) \cdot \nabla U^{n+1} \|_{L^2}
+ C\|g(U^n,\nabla U^n,x,t)\|_{L^2}
\nn\\
& \leq C \| \nabla U^n \|_{L^\infty} \| \nabla
U^{n+1} \|_{L^2}
+ C
 \leq  C,
\label{h2}
\end{align}
for $0\leq n\leq k$.
By the Sobolev interpolation inequality,
\begin{align}
& \|e^{k+1}\|_{L^\infty}
\leq  C\| e^{k+1}\|_{L^2}^{1-d/4}\|
e^{k+1}\|_{H^2}^{d/4}
\leq C\tau^{1-d/4} .
\label{1-infty}
\end{align}

Again we multiply the equation
(\ref{err-TDparab-01}) by $-\Delta e^{n+1}$ to
get
\begin{align}
&D_\tau\biggl(\frac{1}{2} \| \nabla e^{n+1}
\|_{L^2}^2 \biggl)
+\biggl(\sigma(U^n)\Delta e^{n+1}, \, \Delta
e^{n+1} \biggl)
\nn \\
& \le C\epsilon^{-1} \| \nabla \sigma(U^n) \cdot
\nabla e^{n+1} \|_{L^2}^2
+ C \epsilon^{-1} (
\|R^{n+1}\|_{L^2}^2+\|\nabla\cdot([\sigma(u^n)-\sigma(U^n)]
\nabla u^{n+1})\|_{L^2}^2 \nn\\
&~~~+\|g(u^n,\nabla u^n,x,t)-g(U^n,\nabla
U^n,x,t)\|_{L^2}^2\big)
+ \epsilon \| \Delta e^{n+1} \|_{L^2}^2 \, .
\label{fdsjkiolafehio}
\end{align}
By (\ref{g-sigma}), the Sobolev interpolation
inequality and the induction
assumption, we have
\begin{align*}
& \| \nabla \sigma(U^n) \cdot \nabla e^{n+1}
\|_{L^2}
\leq C\| \nabla U^{n} \|_{L^\infty} \| \nabla
e^{n+1} \|_{L^2}
\leq C\|\nabla e^{n+1}\|_{L^2},
\nn \\
& \|g(u^n,\nabla u^n,x,t)-g(U^n,\nabla
U^n,x,t)\|_{L^2} \leq
C \| e^n \|_{H^1}
\end{align*}
and
\begin{align*}
& \| \nabla \cdot \big [
(\sigma(u^{n})-\sigma(U^n))\nabla
u^{n+1}\big]\|_{L^2}
\nn \\
& \leq C\| e^{n+1}\|_{H^1} \| \nabla u^{n+1}
\|_{L^\infty} +
C\| e^{n+1} \|_{L^3} \| \Delta u^{n+1} \|_{L^6}
\nn \\
& \leq C \| e^{n+1}\|_{H^1}
 \, .
\end{align*}
Using Lemma \ref{h23} and choosing a small
$\epsilon$, the inequality
(\ref{fdsjkiolafehio}) reduces to
\begin{align*}
&D_\tau\biggl( \| \nabla e^{n+1} \|_{L^2}^2
\biggl)
+ \| e^{n+1}\|_{H^2}^2
\leq C\|e^{n+1}\|_{H^1}^2 +C\|e^n\|_{H^1}^2+C \|
R^{n+1} \|_{L^2}^2 .
\end{align*}
With Gronwall's inequality, we see that there
exists a positive constant $\tau_2$
such that when $\tau<\tau_2$,
\begin{align}\label{asdfasgghh}
& \|e^{n+1} \|_{H^1}^2
+\sum_{m=0}^{n}\tau\big\|
e^{m+1}\big\|_{H^2}^2\leq C\tau^2
\end{align}
which together with (\ref{regularity}) leads to
\begin{align}
\begin{array}{ll}
&\| e^{n+1}\|_{H^1}\leq C_2\tau,\quad
\|U^{n+1}\big\|_{H^2}
\leq C_2 , \\[5pt]
&\| D_{\tau} U^{n+1} \|_{H^1} \leq C_2,\quad
\sum_{m=0}^{n}\tau\big\|D_\tau
U^{m+1}\big\|_{H^2}^2
\leq C_2,
\end{array}
\label{eh1}
\end{align}
for $0\leq n\leq k$.

Moreover, we rewrite the equation
(\ref{e-TDparab-1}) as
\begin{align}
-\Delta U^{n+1} = \frac{1}{\sigma(U^n)} \left (
g(U^n,\nabla U^n,x,t_n) - D_{\tau} U^{n+1} +
\nabla \sigma(U^n)\cdot
\nabla U^{n+1} \right ) \, .
\label{UU}
\end{align}
By Lemma \ref{h23} and \refe{g-sigma},
\begin{align}\label{sdfjkl23135}
\|U^{n+1}\|_{H^3}
& \leq
\Big\| \frac{1}{\sigma(U^n)} \left ( g(U^n,\nabla
U^n,x,t_n)
- D_{\tau} U^{n+1} + \nabla \sigma(U^n)\cdot
\nabla U^{n+1} \right ) \Big\|_{H^1}
\nn \\
& \leq C\|U^n\|_{H^1}+
C\|g(U^n,\nabla U^n,x,t_n)\|_{H^1} + \|D_\tau
U^{n+1}\|_{H^1}
+C\|\nabla \sigma(U^n) \cdot \nabla
U^{n+1}\|_{H^1}
\nn \\
& \le
C+C\| U^n \|_{H^2}+C\| U^n \|_{H^2}\|\nabla
U^{n+1} \|_{L^\infty}
+ C\|\nabla U^n \|_{L^\infty} \| U^{n+1}\|_{H^2}
\nn \\
& \leq
C +C \|\nabla  U^{n+1} \|_{L^\infty}
\nn \\
& \leq C+\epsilon \|U^{n+1}\|_{H^3}
+C\epsilon^{-1}\|U^{n+1}\|_{H^2} ,
\end{align}
which in turn implies that
\begin{align}
\|U^{n+1}\|_{H^3} \le C_3
\label{Uh3}
\end{align}
if we choose $\epsilon \leq 1/2$.
By the Sobolev interpolation inequality,
\begin{align}
& \|\nabla e^{k+1}\|_{L^\infty}
\leq  C\| e^{k+1}\|_{H^1}^{1-d/4}\|
e^{k+1}\|_{H^3}^{d/4}
\leq C\tau^{1-d/4}
\nn
\end{align}
which with \refe{1-infty} shows that
there exists $\tau_3>0$ such that
\begin{align}
\| U^{k+1} \|_{L^\infty} + \| \nabla U^{k+1}
\|_{L^\infty}
\le \| u^{k+1} \|_{L^\infty} + \| \nabla u^{k+1}
\|_{L^\infty}
+ \| e^{k+1} \|_{L^\infty} + \| \nabla e^{k+1}
\|_{L^\infty} \le K
\,
\end{align}
for $\tau<\tau_3$. Thus (\ref{fsdioo123}) holds
for $n=k+1$ when
$\tau < \tau_0 := \min \{ \tau_1, \tau_2,
\tau_3\}$ and the induction is closed.
  From \refe{eh1}-\refe{Uh3}, we see that
(\ref{ErrestDisSol211})-(\ref{TDErrEstLemm}) hold
when $\tau < \tau_0$.

Secondly, we prove that
(\ref{ErrestDisSol211})-(\ref{TDErrEstLemm}) hold
for
$\tau\geq\tau_0$.
We assume that $\max_{1\leq n\leq k}\|U^n\|_{H^3}\leq \gamma_k$ for
some positive constant
$\gamma_k$ (which may depend upon $\tau_0$)
since $\| U^0 \|_{H^3} = \| u_0 \|_{H^3} \leq C$.
From \refe{e-TDparab-1}, it is easy to see that
\begin{align*}
\| U^{k+1} \|_{L^2} \leq \biggl(\sum_{n=1}^k C \| g(U^n,\nabla
U^n,x,t_n) \|_{L^2}^2\tau\biggl)^{\frac{1}{2}} \leq C_{\gamma_k}
\, .
\end{align*}
Then we
apply Lemma \ref{h23} to \refe{UU}. Via a similar approach as (\ref{sdfjkl23135}), we can
derive that
\begin{align*}
\|U^{n+1}\|_{H^3}\leq C_{\gamma_k}: =\gamma_{k+1} .
\end{align*}
Since $N=T/\tau\leq T/\tau_0$, we take
$C_4=\max_{0\leq k\leq N} \gamma_k$ so that
\begin{align*}
\|U^{n+1}\|_{H^3}\leq  C_4 ,\quad 0\leq n\leq N-1
\end{align*}
which further shows that
\begin{align*}
& \| D_{\tau} U^{n+1} \|_{H^3} \le C_4 \tau_0^{-1}
,\\
&\|u^{n+1}-U^{n+1}\|_{H^1}\leq
\|u^{n+1}\|_{H^1}+\|U^{n+1}\|_{H^1}
\leq M + C_4 .
\end{align*}
Thus the induction is complete for $\tau \geq
\tau_0$.

Combining the two cases, we complete the proof of
Theorem \ref{ErrestDisSol}.
\endproof

\subsection{The proof of Theorem
\ref{ErrestFEMSol} }
For $n\geq 0$, let $R_h^{n+1}: H_0^1(\Omega)
\rightarrow V_h$
be a projection defined by
 \begin{align}\label{projection}
&\big(\sigma(U^n)\nabla ( w- R_h^{n+1}w),\nabla
v\big) =0
\end{align}
for all $v\in V_h$, and we set
$R^0_h:=R^1_h$.
With the regularity of $U^n$ proved in Theorem
\ref{ErrestDisSol},
we have the following inequalities:
\begin{align}
&\|U^{n+1}-R_h^{n+1}U^{n+1}\|_{L^2}
+h\|U^{n+1}-R_h^{n+1}U^{n+1}\|_{H^1}\leq
C\|U^{n+1}\|_{H^3}h^3 , ~~~ \mbox{if}~~ r\geq 2
\label{proj-01} \\
&\|U^{n+1}-R_h^{n+1}U^{n+1}\|_{L^6}
+h\|U^{n+1}-R_h^{n+1}U^{n+1}\|_{W^{1,6}}\leq
C\|U^{n+1}\|_{W^{2,6}}h^2 ,~~~\mbox{if}~~ r=1
 \label{proj-02} \\
&\| R_h^{n+1}U^{n+1}\|_{W^{1,\infty}}\leq
C\|U^{n+1}\|_{W^{1,\infty}} ,
\label{proj-03} \\
& \| D_{\tau} R_h^{n+1} U^{n+1} \|_{W^{1,6}} \le C
\| D_{\tau} U^{n+1} \|_{W^{1,6}},
\label{proj-04} \\
&\| D_{\tau} (U^{n+1}-R_h^{n+1}U^{n+1})
\|_{H^{-1}}
\leq Ch^3 , \quad\mbox{for}~r\geq 2,
\label{proj-05} \\
&\| D_{\tau} (U^{n+1}-R_h^{n+1}U^{n+1}) \|_{L^2}
\leq Ch^2 , \quad~~\mbox{for}~r=1,
\label{proj-06}
\end{align}
where (\ref{proj-01})-(\ref{proj-02}) are standard
error estimates of
elliptic equations,
(\ref{proj-03})-(\ref{proj-04}) follow from
\cite{RS}
and the references therein,
(\ref{proj-05})-(\ref{proj-06}) can be proved
in a similar way as in \cite{LS1,LS2}.

The following inverse inequalities will also be
used in our proof.
\begin{align*}
& \| v \|_{L^p} \le C h^{\frac{d}{p} -
\frac{d}{q}} \| v \|_{L^q},
\quad v \in V_h, ~ 1\leq q\leq p\leq \infty ,
\\
& \| \nabla v \|_{L^p} \le C h^{-1} \| v \|_{L^p},
\quad v \in V_h, ~ 1\leq p\leq \infty .
\end{align*}

Let
$$
K_1=\max_{0\leq n\leq N}\|U^n\|_{W^{1,\infty}}
+\max_{0\leq n\leq N}\| R_h^n
U^n\|_{W^{1,\infty}}+2, \qquad
r^*=\min \{ r, 2\} \, .
$$
By the regularity assumptions for $\sigma$ and
$g$, there exist
$\sigma_{K_1}^*$ and $C_{K_1}^*>0$ such that
\begin{align}
&|\sigma(s)|+|\sigma'(s)|+|\sigma''(s)|+
|g(\alpha,\beta,x,t)|
+|\partial_{\alpha}g(\alpha,\beta,x,t)|
\nn \\
&
+\mbox{$\sum_{j=1}^d$}|\partial_{\beta_j}g(\alpha,\beta,x,t)|
+\mbox{$\sum_{j=1}^d$}
|\partial_{\alpha}\partial_{\beta_j}g(\alpha,\beta,x,t)|
\nn \\
&+\mbox{$\sum_{j=1}^d$}|\partial_{x_j}
\partial_{\beta_j}g(\alpha,\beta,x,t)|
+\mbox{$\sum_{i,j=1}^d$}
|\partial_{\beta_i}\partial_{\beta_j}
g(\alpha,\beta,x,t)|
\leq C_{K_1}^*, \label{K122}\\[5pt]
&\sigma(s)\geq \sigma_{K_1}^* ,\label{K1asfd}
\end{align}
for all $s\in[-K_1,K_1]$ and
$(\alpha,\beta,x,t)\in Q_{K_1}$.

\bigskip

\noindent{\it Proof of Theorem
\ref{ErrestFEMSol}}. ~~~
We shall prove
\begin{align}
&\|U^n_h\|_{L^\infty}+\|\nabla U^n_h\|_{L^\infty}
<K_1, \quad n=0,1,...,N,\label{Uh-bound}\\
&\| e_h^{n} \|_{L^2}^2 + \sum_{m=0}^n \tau \|
\nabla e_h^{m} \|_{L^2}^2
\le \widehat C_0 h^{2r^*+2} ,
\label{r>1}
\end{align}
simultaneously by mathematical induction, where
$e_h^{n}=R_h^{n}U^{n}-U^{n}_h$.
It is easy to see that
the inequalities \refe{Uh-bound}-\refe{r>1} hold
for $n=0$.
So we can assume that \refe{Uh-bound}-(\ref{r>1})
hold for $0\leq n\leq k$.
By (\ref{K122})-(\ref{K1asfd}), the coefficient
matrix of the linear system
(\ref{FDFEM}) is symmetric and positive definite.
Therefore, the system (\ref{FDFEM}) admits a
unique solution $U^{n+1}_h$ in $V_h$
for $0\leq n\leq N-1$.

Since the solution of the time-discrete equation
(\ref{e-TDparab-1}) $U^n$ satisfies
\begin{align}
&\big(D_\tau U^{n+1},v\big)+\big(\sigma(U^n)\nabla
R^{n+1}_hU^{n+1},\nabla v\big)
=\big(g(U^n,\nabla U^n,x,t),v\big)  ,\quad
\forall~v\in V_h ,
\nn
\end{align}
it follows that $e_h^{n+1}$ satisfies the
equation
\begin{align}
&\big(D_\tau
e^{n+1}_h,v\big)+\big(\sigma(U^n_h)\nabla
e^{n+1}_h,\nabla v\big)
\label{errFEMFDFEM}\\
&=-\big((\sigma(U^n)-\sigma(U^n_h))\nabla
R_h^{n+1}U^{n+1},\nabla v\big)
+ \big(D_\tau (R^{n+1}_hU^{n+1}-U^{n+1}),v\big)
\nn \\
& +\big(g(U^n,\nabla U^n,x,t)-g(U^n_h,\nabla
U^n_h,x,t),v\big) ,
\qquad\forall~ v\in V_h.
\nn
\end{align}
Now we estimate the last three terms in the above
equation,
respectively. For the first two terms, we see
that
$$
|\big((\sigma(U^n)-\sigma(U^n_h))\nabla
R_h^{n+1}U^{n+1},\nabla v\big) |
\le C(\| e_h^n \|_{L^2} + h^{r^*+1} ) \| \nabla v
\|_{L^2}
$$
and
$$
| \big(D_\tau (R^{n+1}_hU^{n+1}-U^{n+1}),v\big) |
\le \| D_\tau (R^{n+1}_hU^{n+1}-U^{n+1})
\|_{H^{-1}} \|v \|_{H^1}
\, .
$$
We rewrite the third term by
\begin{align}
&\big(g(U^n,\nabla U^n,x,t)-g(U^n_h,\nabla
U^n_h,x,t),v\big)
\label{g-1}
\\
& = \big(g(U^n,\nabla U^n_h,x,t)-g(U^n_h,\nabla
U^n_h,x,t),v\big)
+\big(g(U^n,\nabla U^n,x,t)-g(U^n,\nabla
U^n_h,x,t),v\big) \, .
\nn
\end{align}
Here, we have
$$
\big|\big(g(U^n,\nabla U^n_h,x,t)-g(U^n_h,\nabla
U^n_h,x,t),v\big)\big|
\le C(\| e_h^n \|_{L^2} + h^{r^*+1}) \| v
\|_{L^2}
\,.
$$
By Taylor's formula, we get
\begin{align*}
& g(U^n,\nabla U^n,x,t)-g(U^n,\nabla U^n_h,x,t)
\nn \\
& = \nabla (U^n-U^n_h) \cdot
\int_0^1 \nabla_2 g(U^n,(1-s)\nabla U^n+s\nabla
U^n_h,x,t)ds
\nn \\
& = \nabla (U^n-U^n_h) \cdot
\int_0^1 \left [ \nabla_2 g(U^n,(1-s)\nabla
U^n+s\nabla U^n_h,x,t)
- \nabla_2 g(U^n,\nabla U^n,x,t) \right ]ds
\nn \\
&~~~ + \nabla\cdot \left ( (U^n-U^n_h)
 \nabla_2 g(U^n,\nabla U^n,x,t)  \right )
 -  (U^n-U^n_h)
\nabla \cdot\left (  \nabla_2 g(U^n,\nabla
U^n,x,t)  \right )
\end{align*}
where $\nabla_2 g$ denotes the gradient of $g$
with respect to the
second conponent. Therefore,
\begin{align}
&\Big | \Big(g(U^n,\nabla U^n,x,t)-g(U^n,\nabla
U^n_h,x,t),v\Big) \Big |
 \nn \\
&~~~~ \le \| \nabla (U^n-U^n_h) \|^2_{L^{12/5}} \|
v \|_{L^6}
+ \| U^n-U^n_h \|_{L^2} \| v \|_{L^2}
\nn \\
&~~~~~~+ \Big | \Big ( (U^n-U^n_h)
 \nabla_2 g(U^n,\nabla U^n,x,t) \, , \nabla v \Big
) \Big |
\nn \\
&~~~~ \leq C (\| \nabla e^n_h \|^2_{L^{12/5}} +
h^{r^*+1}) \| v \|_{H^1}
+ C (\| e^n_h \|_{L^2}  + h^{r^*+1})\| v \|_{H^1}
\, .
\nn
\end{align}
Substituting $v=e^{n+1}_h$ into
(\ref{errFEMFDFEM}),
we obtain
\begin{align}
&D_\tau\Big(\frac{1}{2}\|e^{n+1}_h\|_{L^2}^2\Big)
+ \sigma_{K_1}\|\nabla e^{n+1}_h\|_{L^2}^2
\label{error-1}
\\
& \leq
\epsilon \|\nabla e^{n+1}_h\|_{L^2}^2 + \epsilon
\|e^{n+1}_h\|_{L^2}^2
+ C_6 \epsilon^{-1} \| e^n_h \|_{L^2} + C_7
\epsilon^{-1}\| \nabla e^n_h \|^4_{L^{12/5}}
\nn \\
&~~~+ C\epsilon^{-1} \|D_\tau
(R^{n+1}_hU^{n+1}-U^{n+1})\|_{H^{-1}}^2
+ C\epsilon^{-1} h^{2r^*+2} \, .
\nn
\end{align}
By an inverse inequality and the induction
assumption,
\begin{align}\label{ajlfh9oo88}
\|\nabla e^n_h \|_{L^{12/5}}^4
\leq C \widehat C_0 h^{4-d/3} \| \nabla e_h^n
\|_{L^2}^2
\leq h \| \nabla e_h^n \|_{L^2}^2
\end{align}
if $C\widehat C_0 h^{3-d/3}<1$.
With this estimate, we take
$h \le \epsilon \sigma_{K_1}/(8C_7)$ and sum up
\refe{error-1} to get
\begin{align}
\frac{1}{4}\|e^{n+1}_h\|_{L^2}^2
+ \sum_{m=0}^n \frac{\tau\sigma_{K_1}}{4}\|\nabla
e^{m+1}_h\|_{L^2}^2
\leq C_6 \epsilon^{-1} \sum_{m=0}^{n-1}
\|e^{m+1}_h\|_{L^2}^2
+ C  h^{2r^*+2}
\nn
\end{align}
By (explicit) Gronwall's inequality, we derive
that
\begin{align}\label{error-l2}
\|e^{n+1}_h\|_{L^2} ^2+ \sum_{m=0}^n\tau \|\nabla
e^{m+1}_h\|_{L^2}^2
\leq  C_8 h^{2r^*+2}
\end{align}
for $0\leq n\leq k$.

To complete the mathematical induction, we need to
prove \refe{Uh-bound}-\refe{r>1} for $n=k+1$. For
this purpose, we consider two cases.

{\it Case I:} $r\geq 2$.~
In this case, $r^*=2$ and we can apply inverse
inequalities for \refe{error-l2} to get
\begin{align*}
& \|\nabla e^{n+1}_h\|_{L^\infty}\leq
Ch^{-d/2-1}\|e^{n+1}_h\|_{L^2}\leq C_9h^{r^*-3/2}
,
\\
&\|e^{n+1}_h\|_{L^\infty}\leq
Ch^{-d/2}\|e^{n+1}_h\|_{L^2}\leq C_9h^{r^*-1/2} ,
\end{align*}
which implies that
\begin{align}\label{K1-1}
& \|U^{n+1}_h\|_{W^{1,\infty}}
\leq \|e^{n+1}_h\|_{W^{1,\infty}}
+\|R_h^{n+1}U^{n+1}_h\|_{W^{1,\infty}}  \leq K_1
\end{align}
when $C_9 h^{1/2} < 1$. This completes the
induction for $r\geq 2$, and
\refe{Uh-bound}-\refe{r>1} hold for all $0\le n
\le N-1$ and $0<\tau \le T$.

{\it Case II:} $r=1$.~
In this case, $r^*=1$. To get the boundedness of
$\| \nabla e_h^{n+1} \|_{L^\infty}$, we present
the $H^1$-estimate with
an additional induction assumption:
\begin{equation}
\|\nabla e_h^n \|_{L^2} + \| D_{\tau} e_h^n
\|_{L^2}
\le \widehat C_1 h^2
\, .
\label{r=1}
\end{equation}
From the initial condition, we see that
\refe{r=1} holds for $n=0$ and we can assume that
it holds for $0\leq n\le k$.
We substitute $v=D_{\tau} e^{n+1}_h$ into
(\ref{errFEMFDFEM}). With a similar approach
to \refe{g-1}, we obtain
\begin{align}
&D_\tau\biggl(\Big\|\sqrt{\sigma(U_h^n)}\nabla
e^{n+1}_h\Big\|_{L^2}^2\biggl)
+ \|D_{\tau} e^{n+1}_h\|_{L^2}^2
\label{e-1-1}
\\
&\leq
\big ( D_{\tau} \sigma(U_h^n) \nabla e_h^n, \nabla
e_h^n \big )
- \big((\sigma(U^n)-\sigma(U^n_h))\nabla
R_h^{n+1}U^{n+1}, D_{\tau} \nabla e_h^{n+1}\big)
\nn \\
&~~+ C ( \|D_\tau (R^{n+1}_hU^{n+1}-U^{n+1})
\|_{L^2}
+ \| \nabla e_h^n \|_{L^4}^2  + \| e_h^n \|_{L^2}
+ h^2)
\| D_{\tau} e_h^{n+1} \|_{L^2}
\nn \\
&~~ + \big ( (U^n-U^n_h)   \nabla_2 g(U^n,\nabla
U^n,x,t)  \, ,
D_{\tau} \nabla e_h^{n+1} \big ) \,  .
\nn
\end{align}
Using (\ref{proj-03}) and the Sobolev embedding
inequalities
$$
\|D_\tau R^n_hU^n\|_{L^\infty}\leq C\|D_\tau
R^n_hU^n\|_{W^{1,6}},
~~~~~
\|D_\tau U^n\|_{W^{1,6}}\leq C\|D_\tau
U^n\|_{H^2} .
$$

The first two terms of the right-hand side of the
equation \refe{e-1-1} are bounded by
\begin{align}
\big ( D_{\tau} \sigma(U_h^n) \nabla e_h^n, \nabla
e_h^n \big )
& \leq \| D_{\tau} U_h^n \|_{L^\infty} \| \nabla
e_h^n \|_{L^2}^2
\nn \\
&\le (\| D_{\tau} e_h^n \|_{L^\infty} + \|
D_{\tau} R_h^n U^n \|_{L^\infty})
\| \nabla e_h^n \|_{L^2}^2
\nn \\
& \le C( h^{-3/2}\| D_{\tau} e_h^n \|_{L^2} +  \|
D_{\tau}  U^n \|_{H^2}
) \| \nabla e_h^n \|_{L^2}^2
\nn
\end{align}
and
\begin{align*}
& -\big((\sigma(U^n)-\sigma(U^n_h))\nabla
R_h^{n+1}U^{n+1},
D_{\tau} \nabla e_h^{n+1}\big) \nn \\
&~~~~
= -D_{\tau} \big((\sigma(U^n)-\sigma(U^n_h))\nabla
R_h^{n+1}U^{n+1}, \nabla e_h^{n+1}\big)
\nn \\
&~~~~~~
+\big(D_{\tau}[(\sigma(U^n)-\sigma(U^n_h))\nabla
R_h^{n+1}U^{n+1}], \nabla e_h^n\big)
\nn \\
&~~~~ \le
-D_{\tau} \big((\sigma(U^n)-\sigma(U^n_h))\nabla
R_h^{n+1}U^{n+1}, \nabla e_h^{n+1}\big)
\nn \\
&~~~~~~
+ C(\|D_{\tau}(U^n-U^n_h) \|_{L^2} +
\|(U^n-U^n_h)D_{\tau}U^n \|_{L^2})
\| \nabla R_h^nU^{n} \|_{L^\infty}\|\nabla
e_h^{n+1} \|_{L^2}
\nn \\
&~~~~~~ + C \|U^n-U^n_h \|_{L^6} \| D_{\tau}
\nabla R_h^{n+1}U^{n+1}] \|_{L^3}
\|\nabla e_h^{n+1} \|_{L^2}
\nn \\
&~~~~ \leq
- D_{\tau} \big((\sigma(U^n)-\sigma(U^n_h))\nabla
R_h^{n+1}U^{n+1}, \nabla e_h^{n+1}\big)
\nn \\
&~~~~~~ + \epsilon \|D_{\tau} e_h^n \|_{L^2}^2 +
C\epsilon^{-1} (\| D_{\tau} U^{n+1} \|_{H^2}^2 +1
)
\|\nabla e_h^{n+1} \|_{L^2}^2 +
\epsilon\|e^n_h\|_{H^1}^2 + Ch^4 .
\end{align*}
Moreover, we have
\begin{align*}
& \big ( (U^n-U^n_h)
\nabla_2 g(U^n,\nabla U^n,x,t)  \, , D_{\tau}
\nabla e_h^{n+1} \big )
\nn \\
& = D_{\tau} \big ( (U^n-U^n)
 \nabla_2 g(U^n,\nabla U^n,x,t) \, , \nabla
e_h^{n+1} \big )\\
&~~~ - \big ( D_{\tau} [(U^n-U^n_h)
 \nabla_2 g(U^n,\nabla U^n,x,t)] \, , \nabla e_h^n
\big )
\nn \\
& \le
D_{\tau} \big ( (U^n-U^n_h)
 \nabla_2 g(U^n,\nabla U^n,x,t) \, , \nabla
e_h^{n+1} \big )
+ \epsilon (\| D_{\tau} e^n_h \|_{L^2}^2+ \|
D_{\tau} (U^n-R^n_hU^n) \|_{L^2}^2)
\nn \\
&~~~ + C\epsilon^{-1}(1+\|D_\tau U^n\|_{H^2}) \|
\nabla e_h^n \|_{L^2}^2
+ \epsilon (\| e_h^n \|_{L^6}^2+ \|U^n-U_h^n
\|_{L^6}^2) \, .
\end{align*}
With the above estimates, the inequality
(\ref{e-1-1}) reduces to
\begin{align}\label{r1-1}
&D_\tau\biggl(\Big\|\sqrt{\sigma(U_h^n)}\nabla
e^{n+1}_h\Big\|_{L^2}^2\biggl)
+\|D_{\tau} e^{n+1}_h\|_{L^2}^2
\\
&\leq C \epsilon^{-1} (1 + \| D_{\tau} U^n
\|_{H^2}^2)
(\| \nabla e_h^n \|_{L^2}^2 + \| \nabla e_h^{n+1}
\|_{L^2}^2)
+ \epsilon (\|D_{\tau} e_h^n \|_{L^2}^2+\|e_h^n
\|_{H^1}^2)
\nn \\
&~~~ +  C\epsilon^{-1} \|D_\tau
(R^{n+1}_hU^{n+1}-U^{n+1})
\|_{L^2}^2+C\epsilon^{-1}\|D_\tau (R^n_hU^n-U^n)
\|_{L^2}^2
 \nn\\
&~~~+ D_{\tau} J^{n+1} + C h^4
\nn
,
\end{align}
where
\begin{align*}
J^{n+1} =
& -\big((\sigma(U^n)-\sigma(U^n_h))\nabla
R_h^{n+1}U^{n+1}, \nabla e_h^{n+1}\big)
\nn \\
&+ \big ( (U^n-U^n)
 \nabla_2 g(U^n,\nabla U^n,x,t) \, , \nabla
e_h^{n+1} \big )
\, .
\end{align*}
Furthermore, summing up \refe{r1-1} gives
\begin{align}\label{r1-2}
& \| \nabla e^{n+1}_h\|_{L^2}^2 + \sum_{m=0}^n
\tau \| D_{\tau} e_h^{m+1} \|_{L^2}^2 \\
&\leq
C_{11}  \sum_{m=0}^n \tau
\left (1+ \| D_{\tau} U^{m+1} \|_{H^2}^2 \right)
\| \nabla e^{m+1}_h\|_{L^2}^2
+ C h^4
\nn
\end{align}
where we have noted that
\begin{align*}
& | J^{n+1} | \leq \epsilon \| \nabla
e^{n+1}_h\|_{L^2}^2 + C\epsilon^{-1} h^4 .
\end{align*}
Since
\begin{align*}
& \sum_{m=0}^n \tau
\| D_{\tau} U^{m+1} \|_{H^2}^2
\leq C,
\end{align*}
by Gronwall's inequality and
(\ref{ErrestDisSol211}),
(\ref{r1-2}) further reduces to
\begin{align*}
\| \nabla e^{n+1}_h\|_{L^2}^2
+ \sum_{m=0}^n \tau \| D_{\tau} e_h^{m+1}
\|_{L^2}^2
\leq C_{12} h^4
\end{align*}
for $0\leq n\leq k$, provided $\tau <\tau_5$ for
some positive constant $\tau_5$.

Now by an inverse inequality, we have the
estimate
$$
\| e^{n+1}_h\|_{W^{1,\infty}} \leq C h^{-d/2}
\|e^{n+1}_h\|_{H^1} \leq C_{13} h^{1/2},
$$
and so
\begin{align}\label{K1}
& \|U^{n+1}_h\|_{W^{1,\infty}}
\leq  \|e^{n+1}_h\|_{W^{1,\infty}} +
\|R_h^{n+1}U^{n+1}\|_{W^{1,\infty}}  \leq K_1
\nn
\end{align}
when $C_{13}h^{1/2}<1$. It suffices to choose
$\widehat C_0\geq 1+C_8$
and $\widehat C_1\geq 1+2\sqrt{C_{12}}$ so that
the mathematical induction is closed when $\tau <
\tau_5$. It follows that
(\ref{Uh-bound}), (\ref{r>1}) and (\ref{r=1}) hold
for all $n=1,2,...,N$ when
$\tau < \tau_5$.

When $\tau \geq \tau_5$ and $r=1$, we can see from
\refe{error-l2} that
\begin{align}
\| e_h^{n+1} \|_{H^1} \le C_{14}\tau^{-1}_5 h^2
\end{align}
which together with an inverse inequality implies
that
$$
\|e_h^{n+1} \|_{W^{1,\infty}}
\leq Ch^{-d/2}\| e_h^{n+1} \|_{H^1} \leq  C_{15}
\tau_5^{-1} h^{1/2}<1
$$
if $h<C_{15}^2/\tau_5^2$. Therefore,
\begin{align*}
& \| U_h^{n+1} \|_{L^\infty}
+ \| \nabla U_h^{n+1} \|_{L^\infty}
\\
& \qquad
\leq \| R_h^{n+1} U^{n+1} \|_{L^\infty}+  \|
\nabla R_h^{n+1} U^{n+1} \|_{L^\infty}+
\| e_h^{n+1}  \|_{L^\infty} + \| \nabla e_h^{n+1}
\|_{L^\infty} \leq K_1 .
\end{align*}
It suffices to choose $\widehat C_0\geq 1+C_8$ and
$\widehat C_1\geq 1+C_{14}\tau^{-1}_5$
so that the mathematical induction is closed for
$\tau \geq \tau_5$.
It follows that
(\ref{Uh-bound}), (\ref{r>1}) and (\ref{r=1})
hold for all $n=1,2,...,N$ when
$\tau \geq\tau_5$.

Combining the two cases, we complete the proof of
Theorem \ref{ErrestFEMSol}.
\endproof

\vskip0.1in

\noindent{\bf Remark 3.1}~~ We have proved Theorem
\ref{ErrestFEMSol} for any $r$-order
Galerkin FEMs under the regularity
assumption \refe{regularity}. Based on the
classical theory of finite element
approximation and interpolation, this assumption
is enough to obtain
optimal error estimates for linear and quadratic
Galerkin FEMs.
In fact, the optimal $L^2$ error bounds for $r=1$
and $r=2$ have been given in \refe{r>1}.
Since the estimates in \refe{r>1} are
$\tau$-independent, by an inverse inequality,
$$
\| e^{n+1} \|_{H^1} \le C h^{r}
\, .
$$
By Theorem \ref{ErrestDisSol} and the projection
error estimates in \refe{proj-01},
we have optimal error estimates
for the linear and quadratic Galerkin FEMs, which
are summarized below.

\begin{corollary}\label{r=2}
{\it Under the assumptions of Theorem
\ref{ErrestFEMSol}, there exist positive
constants $C$ and $h_0$ such that when $h<h_0$,
\begin{align}
&\| U_h^n - u^n \|_{L^2} \leq C(\tau+ h^{r+1}) \\
&
\| U_h^n - u^n \|_{H^1}\leq C(\tau+ h^{r})
\label{error-r=2}
\end{align}
for $r=1$ or $r=2$.
}
\end{corollary}

\section{Error analysis}
\setcounter{equation}{0}
Based on the boundedness of the numerical solution
proved in the last section,
one can easily obtain optimal error estimates of
any $r$-order Galerkin FEMs under corresponding
regularity
assumptions, by following the classical approach
of FEM analysis.
Also it is possible to present the optimal error
estimate for $e_h^n$
as we did in Section 3.2.
However, this requires a rigorous analysis for
stronger regularity of
the time-discrete system.
For simplicity, we follow the classical FEM
approach and give a brief
proof of optimal error estimates of the fully
discrete Galerkin FEM.
In this section,
we assume that the solution to the
initial-boundary
value problem (\ref{e-parab-1})-(\ref{BC}) exists,
satisfying (\ref{regularity}) and the following
condition
\begin{align}
\label{regularity-2}
&\|u\|_{L^\infty((0,T);H^{r+1})}
+\|\partial_t u\|_{L^2((0,T);H^{\max(r,2)})}
 \leq C .
\end{align}

Let $\theta_h^n = U_h^n - \overline R_h^n u^n$
where $\overline R_h^n$
is the elliptic projection defined by
\begin{align}\label{projection}
&\big(\sigma(u^n)\nabla ( w- \overline
R_h^{n+1}w),\nabla v\big) =0
\end{align}
for all $v\in V_h$, and we set $\overline
R^0_h:=\overline R^1_h$. Easy to see that
(\ref{proj-01})-(\ref{proj-06}) also hold for the
projection operator $\overline R_h^{n+1}$.

\begin{theorem}\label{asdfa123}
{\it Suppose that the problem
\refe{e-parab-1}-\refe{IniC} has a unique
solution $u$ satisfying (\ref{regularity}) and
(\ref{regularity-2}) for some positive integer
$r$. Then
there exists $h_6>0$ such that for $h<h_6$,
\begin{align}
\| U_h^n - u^n \|_{L^2} \le C(\tau+ h^{r+1})
\label{error-r>2}
.
\end{align}
}
\end{theorem}

\noindent{\it Proof}~~~
Note that the error function $\theta_h^n$
satisfies
the following equation:
\begin{align}
&\big(D_\tau
\theta^{n+1}_h,v\big)+\big(\sigma(U^n_h)\nabla
\theta^{n+1}_h,\nabla v\big)
\label{e-1}\\
&=-\big((\sigma(u^n)-\sigma(U^n_h))\nabla
\overline R_h^{n+1}u^{n+1},\nabla v\big)
+ \big(D_\tau (\overline R^{n+1}_h
u^{n+1}-u^{n+1}),v\big)
\nn \\
& ~~~ +\big(g(u^n,\nabla u^n,x,t)-g(U^n_h,\nabla
U^n_h,x,t),v\big) +(R^{n+1},v)~
\quad\mbox{for all $v\in V_h$},
\nn
\end{align}
where $R^{n+1}$ is the truncation error satisfying
(\ref{trdfsjhklh}).

By the same approach as used in the proof of
Theorem \ref{ErrestFEMSol}, we can derive that
\begin{align*}
& |\big((\sigma(u^n)-\sigma(U^n_h))\nabla
\overline R_h^{n+1}u^{n+1},\nabla v\big) |
\le C\| \theta_h^n \|_{L^2} \| \nabla v \|_{L^2},
\\
& | \big(D_\tau (\overline R^{n+1}_h
u^{n+1}-u^{n+1}),v\big) |
\le \| D_\tau (\overline R^{n+1}_h
u^{n+1}-u^{n+1}) \|_{L^2} \|v \|_{L^2},
\nn \\
&\big | \big(g(u^n,\nabla u^n,x,t)-g(U^n_h,\nabla
U^n_h,x,t),v\big) \big |\\
&\le C \| \nabla \theta^n_h \|_{L^2}\| \nabla
\theta^n_h \|_{L^\infty}
\| v \|_{L^2} + Ch^{r+1} \| v \|_{H^1}
+ C (\| \theta^n_h \|_{L^2}  + h^{r+1})\| v
\|_{H^1} \, .
\nn
\end{align*}
Since $\|\nabla\theta^n_h\|_{L^\infty}\leq C$ as
implied by Theorem \ref{ErrestFEMSol}, by taking
$v = \theta_h^{n+1}$ in (\ref{e-1}), we derive
that
\begin{align*}
&D_\tau\Big(\frac{1}{2}\|\theta^{n+1}_h\|_{L^2}^2
\Big)
+ \sigma_{K_1}\|\nabla \theta^{n+1}_h\|_{L^2}^2
\\
&\leq
\epsilon( \|\nabla\theta^{n+1}_h\|_{L^2}^2
+\|\nabla\theta^{n+1}_h \|_{L^2})+ C\epsilon^{-1}
(\|\theta^{n}_h\|_{L^2}^2+ \|
\theta_h^{n+1} \|_{L^2}^2)
\nn \\
&~~~
+ C\epsilon^{-1} \|D_\tau (\overline R^{n+1}_h
u^{n+1}-u^{n+1})\|_{H^{-1}}^2
+ C(\|R^{n+1}\|_{L^2}^2+ h^{2r+2}) \, .
\nn
\end{align*}
By Gronwall's inequality, there exists
a positive constant $\tau_6$ such that when
$\tau<\tau_6$, we have
\begin{align*}
\max_{0\leq n\leq N}\|\theta^{n}_h\|_{L^2}^2
+ \sum_{m=0}^N \tau \|\nabla
\theta^{m}_h\|_{L^2}^2
\leq C_{17} (\tau ^2 + h^{2r+2}) .
\end{align*}
Therefore, (\ref{error-r>2}) holds when
$\tau<\tau_6$.

For $\tau\geq\tau_6$, by Theorem
\ref{ErrestFEMSol},
$\| U_h^n - u^n \|_{L^2} \leq C\leq
C\tau_6^{-1}(\tau+h^{r+1})$.

The proof of Theorem \ref{asdfa123} is complete.
\endproof

\newpage

\section{Numerical examples}

\setcounter{equation}{0}

\noindent{\it Example 5.1}~~
First, we consider an artificial example governed
by the equation
\begin{align}
&\frac{\partial u}{\partial t}-\Delta u=
\sigma(u)|\nabla u|^4+f ,
\label{dfnhsdfui001}
\end{align}
in the domain $\Omega = (0,1)\times(0,1)$ with
$\sigma(u)=1/(1+u)$.
The function $f$ is
chosen corresponding to the exact solution
\begin{align}
&u(x,y,t)=10x(1-x)y(1-y)
\,{\rm sech}(x+y-t)^2
\label{dfnhs001}
\end{align}
which satisfies the homogeneous Dirichlet boundary
condition.

A uniform triangular partition with $M+1$ nodes in
each direction
is used in our computation (with $h=\sqrt{2}/M$).
We solve the system by the
proposed method with a linear FE approximation up
to the time $t=1$.
To illustrate our error estimates,
we take $\tau=h^2$ and we present numerical
results in
Table \ref{linear-L2-1-1}, from which we can see
that the $L^2$ errors
are proportional to $h^2$.
To demenstrate the unconditional convergence, we
take several different
spatial meshes with  $M= 16, 32, 64$ for each
$\tau=0.01,0.025,0.05$ and we present
numerical errors in Table \ref{linear-L2-2-1}.
Based on our
theoretical analysis, in this case,
\begin{align*}
\| U_h^N - u(\cdot ,t_N) \|_{L^2} = O(\tau + h^2)
\qquad
\| U_h^N - u(\cdot ,t_N) \|_{H^1} = O(\tau + h)
\end{align*}
which tend to $O(\tau)$ as $h\rightarrow 0$.
We can observe from Table \ref{linear-L2-2-1} that
for a fixed $\tau$,
numerical errors behave like $O(\tau)$ as
$h/\tau\rightarrow 0$,
which shows that no time step condition is
needed.

\vskip0.1in

\begin{table}[htp]
\vskip-0.2in
\begin{center}
\caption{$L^2$-norm errors of the linear Galerkin
FEM (Example 5.1). }\vskip 0.1in
\label{linear-L2-1-1}
\begin{tabular}{l|l|c|ccc}
\hline
$\displaystyle\tau=h^2$   &  $h$
& $\displaystyle\| U_h^N - u(\cdot ,t_N)
\|_{L^2}$
& $\displaystyle\| U_h^N - u(\cdot ,t_N) \|_{H^1}$
\\
\cline{2-4}
   & 1/8
& 3.861E-02 &1.657E-01 \\
 & 1/16
& 7.285E-02 &3.211E-02 \\
 & 1/32
& 1.720E-03 &7.678E-03  \\
\hline
\multicolumn{2}{c|}{convergence rate}
& 2.08   & 2.06    \\
\hline
\end{tabular}
\end{center}

\begin{center}
\vskip-0.1in
\caption{$L^2$-norm errors of the linear Galerkin
FEM  with refined meshes (Example 5.1).}
\vskip 0.1in
\label{linear-L2-2-1}
\begin{tabular}{l|c|c|ccc}
\hline
$\tau=0.01$  &  $h$
  & $\| U_h^N - u(\cdot ,t_N) \|_{L^2}$
  & $\| U_h^N - u(\cdot ,t_N) \|_{H^1}$     \\
\cline{2-4}
 & 1/16  & 9.591E-03 &4.526E-02\\
 & 1/32  & 5.484E-03 &3.026E-02\\
 & 1/64  & 4.673E-03 &2.793E-02  \\
\hline
\hline
$\tau=0.025$ &  $h$ 
  & $\| U_h^N - u(\cdot ,t_N) \|_{L^2}$
  & $\| U_h^N - u(\cdot ,t_N) \|_{H^1}$     \\
\cline{2-4}
 & 1/16
 & 1.569E-02  & 8.022E-02\\
 & 1/32
 & 1.167E-02  & 6.700E-02 \\
 & 1/64
 & 1.079E-02  & 6.445E-02 \\
\hline
\hline
$\tau=0.05$ & $h$
  & $\| U_h^N - u(\cdot ,t_N) \|_{L^2}$
  & $\| U_h^N - u(\cdot ,t_N) \|_{H^1}$     \\
\cline{2-4}
 & 1/16
 & 2.486E-02  & 1.312E-01  \\
 & 1/32
 & 2.079E-02  & 1.187E-01  \\
 & 1/64
 & 1.984E-02  & 1.159E-01  \\
\hline
\end{tabular}
\end{center}
\end{table}

\noindent{\it Example 5.2}~~
Secondly, we consider the Burger's equation
\begin{align}
&\frac{\partial {\bf u}}{\partial t}+{\bf
u}\cdot\nabla {\bf u}-\Delta {\bf u}= f ,
\label{dfnhsdfui001}
\end{align}
in the unit disk on the plane, with inhomogeneous
boundary condition ${\bf u}=g$
on $\partial\Omega$. The functions $f$ and $g$ are
given corresponding to the exact solution
\begin{align}
{\bf u}(x,y,t)=({\rm sech}(x+y-t)^2, \,  {\rm
cosh}(x+y-t)^2 ) .
\label{dfnhs001}
\end{align}

The mesh generated here consists of $M$
boundary points with $M=32,64,128$, respectively. See 
Figure 1 for the triangulation of the domain.
Numerical errors with fixed $\tau$ and several
different $h$
are presented in Tables \ref{linear-L2-1-2} and
\ref{linear-L2-2-2}. We can see clearly again
from Table \ref{linear-L2-1-2} that the numerical
errors in $L^2$-norm and $H^1$-norm
are proportional to $O(h^2)$ and $O(h)$,
respectively, when $\tau =O(h^2)$
and from Table \ref{linear-L2-2-2} that
numerical errors behave like $O(\tau)$
as $h/\tau \rightarrow 0$. Thus
no time-step condition is needed.
\vskip0.1in
\begin{figure}[htp]
\begin{minipage}[b]{0.40\linewidth}
\centering
\includegraphics[width=\textwidth]{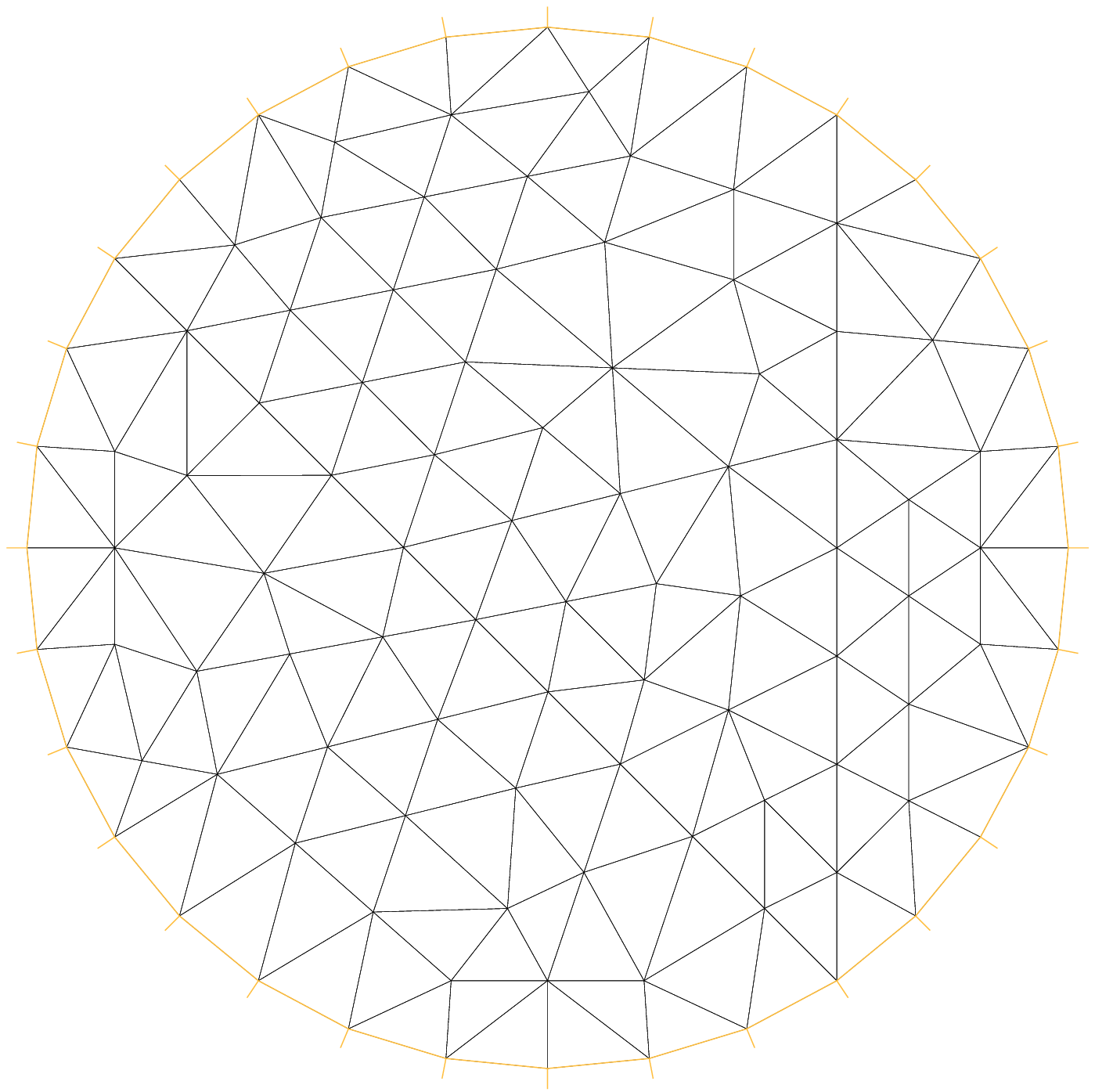}
\end{minipage}
\hspace{-1.8cm}
\begin{minipage}[b]{0.40\linewidth}
\centering
\includegraphics[width=\textwidth]{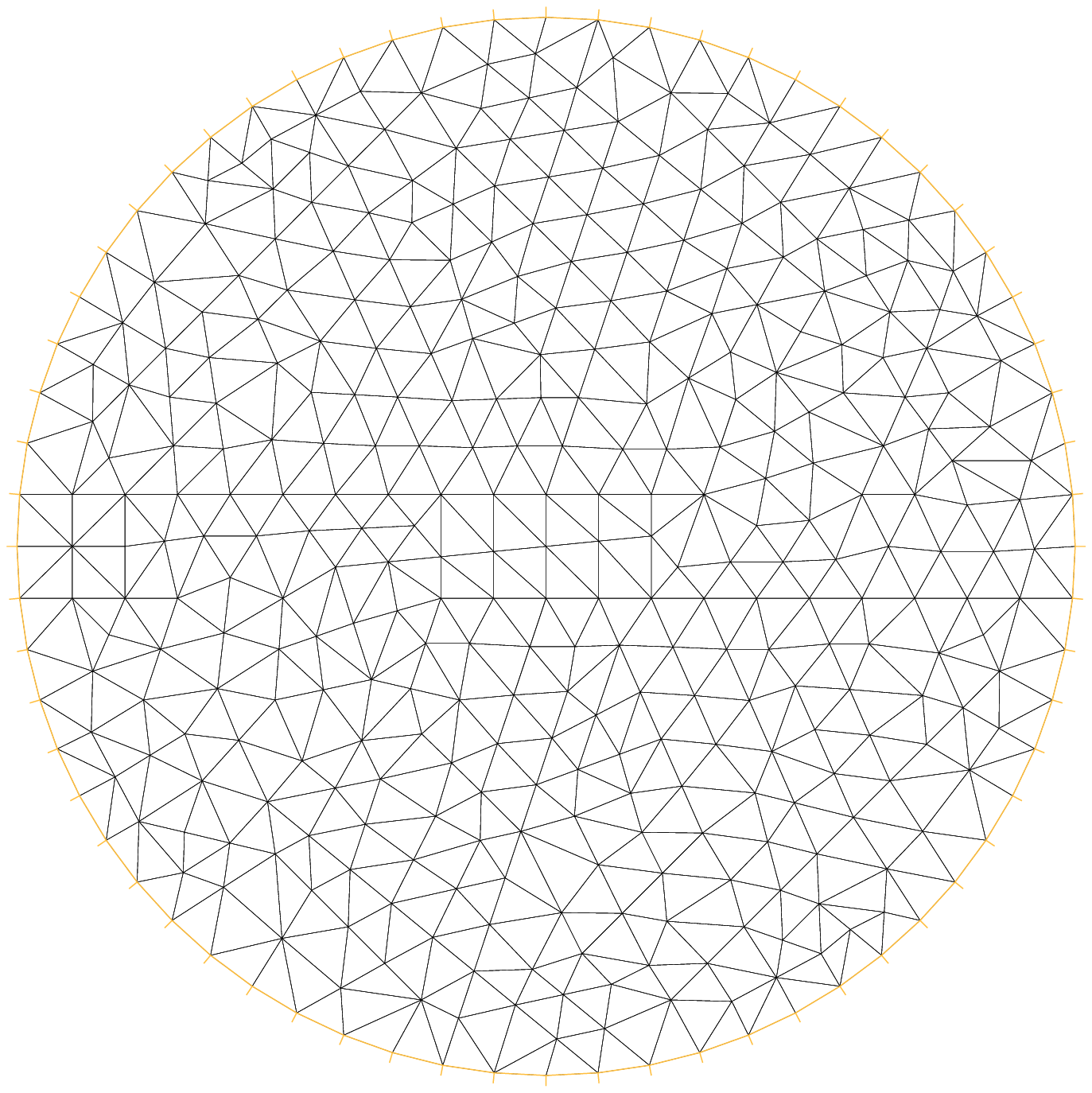}
\end{minipage}
\hspace{-1.8cm}
\begin{minipage}[b]{0.40\linewidth}
\centering
\includegraphics[width=\textwidth]{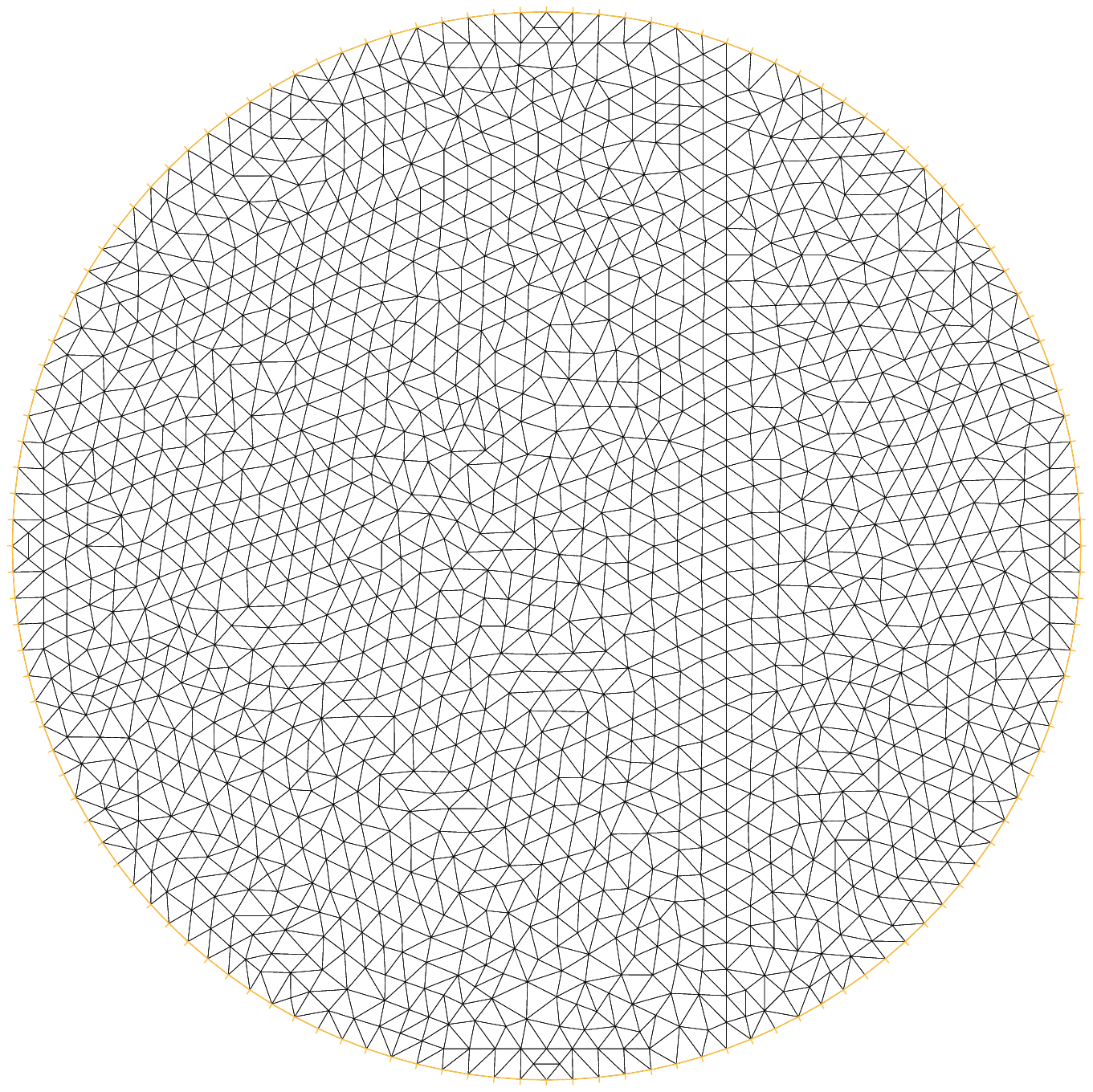}
\end{minipage}
\caption{The FEM meshes with $M=16$, $M=32$ and
$M=64$, respectively.}\label{figure1}
\end{figure}
\begin{table}[t]
\vskip-0.2in
\begin{center}
\caption{$L^2$-norm errors of the linear Galerkin
FEM (Example 5.2).}\vskip 0.1in
\label{linear-L2-1-2}
\begin{tabular}{l|l|c|ccc}
\hline
$\displaystyle\tau=\frac{1}{M^{2}}$  &  $M$
& $\displaystyle\| {\bf U}_h^N - {\bf u}(\cdot
,t_N) \|_{L^2}$
& $\displaystyle\| {\bf U}_h^N - {\bf u}(\cdot
,t_N) \|_{H^1}$    \\
\cline{2-4}
   & 16
& 2.138E-02  & 1.596E-01\\
  & 32
& 4.845E-03  & 7.534E-02\\
 & 64
& 1.314E-03  & 3.633E-02 \\
\hline
\multicolumn{2}{c|}{convergence rate}
& 2.01     &   1.06  \\
\hline
\end{tabular}
\end{center}
\begin{center}
\vskip-0.1in
\caption{$L^2$-norm errors of the linear Galerkin
FEM with  refined meshes (Example 5.2).}
\vskip 0.1in
\label{linear-L2-2-2}
\begin{tabular}{l|c|c|ccc}
\hline
$\tau=0.005$ & $M$
  & $\|{\bf U}_h^N - {\bf u}(\cdot ,t_N)
\|_{L^2}$
  & $\|{\bf U}_h^N - {\bf u}(\cdot ,t_N) \|_{H^1}$
\\
\cline{2-4}
 & 32
 & 2.112E-02  & 1.644E-01 \\
 & 64
 & 5.493E-03  & 8.883E-02 \\
 & 128
 & 4.882E-03  & 5.508E-02\\
\hline
\hline
$\tau=0.010$ &  $M$ 
  & $\|{\bf U}_h^N - {\bf u}(\cdot ,t_N)
\|_{L^2}$
  & $\|{\bf U}_h^N - {\bf u}(\cdot ,t_N) \|_{H^1}$
\\
\cline{2-4}
 & 32
 & 2.078E-02  & 1.897E-01\\
 & 64
 & 9.551E-03  & 1.158E-01 \\
 & 128
 & 1.020E-02  & 8.863E-02 \\
\hline
\hline
$\tau=0.025$  &  $M$
 & $\|{\bf U}_h^N - {\bf u}(\cdot ,t_N)\|_{L^2}$
& $\|{\bf U}_h^N - {\bf u}(\cdot ,t_N)
\|_{H^1}$\\
\cline{2-4}
 & 32   & 2.741E-02 & 2.862E-01\\
 & 64   & 2.520E-02 & 2.228E-01\\
 & 128  & 2.654E-02 & 2.051E-01\\
\hline
\end{tabular}
\end{center}
\end{table}

\noindent{\it Example 5.3}~~
Finally, we consider the equation
\begin{align}
&\frac{\partial u}{\partial t}-\nabla\cdot(\kappa(u)\nabla u)=
\sigma(u)|\nabla u|^4+f
\label{dfnhsdfui00123}
\end{align}
in
$\Omega = (0,1)\times(0,1)\times(0,1)$ with
$\kappa(u)=1+\sin^2u$ and $\sigma(u)=1/(1+u)$.
The function $f$ is chosen corresponding to the
exact solution
\begin{align}
&u(x,y,t)=100x(1-x)y(1-y)z(1-z)
\,\sin(x+2y-z)te^{-t} .
\label{dfnhs001}
\end{align}
which satisfies the homogeneous Dirichlet boundary
condition.

A uniform tetrahedral partition with $M+1$ nodes
in each direction
is used in our computation (with $h=\sqrt{3}/M$).
We solve the system by the
proposed method up to the time $t=1$.
To illustrate our error estimates,
errors of the numerical solution with $\tau=8h^2$
are presented in Table \ref{linear-L2-1}.
Similalry
numerical errors with fixed $\tau$ and refined
$h$
are presented in Table \ref{linear-L2-2}.
The same observations can be made here.
Again, our numerical results show that
the scheme is unconditionally stabe (convergent).
\vskip0.1in
\begin{table}[htp]
\vskip-0.2in
\begin{center}
\caption{$L^2$-norm rrors of the linear Galerkin
FEM (Example 5.3).}\vskip 0.1in
\label{linear-L2-1}
\begin{tabular}{l|l|c|ccc}
\hline
$\tau$   &  $h$
& $\displaystyle\| U_h^N - u(\cdot ,t_N)
\|_{L^2}$
& $\displaystyle\| U_h^N - u(\cdot ,t_N) \|_{H^1}$
\\
\cline{1-4}
1/8   & 1/8
& 2.094E-02 &8.379E-02 \\
1/32   & 1/16
& 4.996E-03 &1.983E-02 \\
1/128  & 1/32
&  1.220E-03 & 4.755E-03  \\
\hline
\multicolumn{2}{c|}{convergence rate}
& 2.08   & 2.06    \\
\hline
\end{tabular}
\end{center}
\begin{center}
\vskip-0.1in
\caption{$L^2$-norm errors of the linear Galerkin
FEM  with refined meshes (Example 5.3).}
\vskip 0.1in
\label{linear-L2-2}
\begin{tabular}{l|c|c|ccc}
\hline
$\tau=0.025$ &  $h$ 
  & $\| U_h^N - u(\cdot ,t_N) \|_{L^2}$
  & $\| U_h^N - u(\cdot ,t_N) \|_{H^1}$     \\
\cline{2-4}
 & 1/8
 & 1.808E-02  & 6.327E-02\\
 & 1/16
 &  4.858E-03  & 1.879E-02 \\
 & 1/32
 &  1.549E-03  & 7.156E-03 \\
\hline
\hline
$\tau=0.05$ & $h$
  & $\| U_h^N - u(\cdot ,t_N) \|_{L^2}$
  & $\| U_h^N - u(\cdot ,t_N) \|_{H^1}$     \\
\cline{2-4}
 & 1/8
 & 1.862E-02  & 6.707E-02  \\
 & 1/16
 & 5.478E-03  & 2.357E-02  \\
 & 1/32
 & 2.181E-03  & 1.206E-02  \\
\hline
\hline
$\tau=0.1$ & $h$
  & $\| U_h^N - u(\cdot ,t_N) \|_{L^2}$
  & $\| U_h^N - u(\cdot ,t_N) \|_{H^1}$     \\
\cline{2-4}
 & 1/8
 & 2.008E-02  & 7.762E-02  \\
 & 1/16
 & 7.241E-03  & 3.759E-02  \\
 & 1/32
 & 3.976E-03  & 2.668E-02  \\
\hline
\end{tabular}
\end{center}
\end{table}

\section{Conclusion}
We have presented unconditionally optimal error
estimates of a class of linearized Galerkin FEMs
for
general nonlinear parabolic equations, which may
cover many physical applications.
The time-step size restriction was always a key
issue in previous analysis and practical
computation. Our theoretical analysis and
numerical results show clearly that no time-step condition is needed for these linearized
Galerkin FEMs.
Our approach is based on a priori estimates of
the numerical solution in the
$W^{1,\infty}$ norm. With these estimates, optimal
error estimates can be proved unconditionally
from classical FEM error analysis. Clearly, our approach
is applicable to many other time discretization schemes
and
more general nonlinear equations (systems).

\end{document}